\documentclass[12pt]{article}
\usepackage{amsfonts}
\usepackage{amssymb}
\usepackage{amsmath}

\def\thebibliograph#1#2{\section*{{\normalsize \bf #2}}\list

 {[\arabic{enumi}]}{\settowidth\labelwidth{[#1]}\leftmargin\labelwidth
     \advance\leftmargin\labelsep
     \usecounter{enumi}}
     \def\newblock{\hskip .11em plus .33em minus -.07em}
     \sloppy
     \sfcode`\.=1000\relax}

\newcommand{\re}{{\mathbb{R}}}
\newcommand{\R}{{\mathbb{R}}^n}
\newcommand{\N}{{\mathbb{N}}}
\newcommand{\Z}{{\mathbb{Z}}}

\newtheorem{thm}{Theorem}
\newtheorem{prop}{Proposition}
\newtheorem{DE}{Definition}

\newtheorem{LEM}{Lemma}
\newtheorem{rem}{Remark}

\newcommand{\be}{\begin{equation}}
\newcommand{\ee}{\end{equation}}
\newcommand{\beq}{\begin{eqnarray}}
\newcommand{\eeq}{\end{eqnarray}}
\newcommand{\beqq}{\begin{eqnarray*}}
\newcommand{\eeqq}{\end{eqnarray*}}

{\vskip\baselineskip\noindent\textbf{Proof of Theorem \protect\ref{#1}\quad}}%

{\vskip\baselineskip\noindent\textbf{Proof of Proposition \protect\ref{#1}\quad}}%

{\vskip\baselineskip\noindent\textbf{Proof of Corollary \protect\ref{#1}\quad}}%

{\vskip\baselineskip\noindent\textbf{Proof of Lemma \protect\ref{#1}\quad}}%

\title{ Continuity of composition operators in Sobolev spaces}
  \author{ G\'erard Bourdaud \& Madani Moussai}
 \date{\today}

\begin{document}

\maketitle

\begin{abstract} We prove that all the composition operators $T_f(g):= f\circ g$, which take the Adams-Frazier space
 $W^{m}_{p}\cap \dot{W}^{1}_{mp}(\re^n)$ to itself, are continuous mappings from $W^{m}_{p}\cap \dot{W}^{1}_{mp}(\re^n)$ to itself, for every integer $m\geq 2$ and every real number $1\leq p<+\infty$. The same automatic continuity property
 holds for Sobolev spaces $W^m_p(\R)$ for $m\geq 2$ and $1\leq p<+\infty$.

\end{abstract}



\section{Introduction}

We want to establish the so-called {\em automatic continuity} property for composition operators in classical Sobolev spaces, i.e. the following statement: 

\begin{thm}\label{contm>0} Let us consider an integer $m>0$, and  $1\leq p<+\infty$.  
If $f:\re \rightarrow \re$ is a function s.t. the composition operator $T_f(g):= f\circ g$
 takes $W^m_p(\R)$ to itself, then $T_f$ is a continuous mapping  from $W^m_p(\R)$ to itself.
 \end{thm}

This theorem has been proved
\begin{itemize} 
\item for $m=1$, see \cite{AA} in case $p=2$, and $ \cite{MM}$ in the general case,
\item for $m>n/p$, $m>1$ and $p>1$ \cite{BL}. 
\end{itemize}
It holds also trivially  in the case of Dahlberg degeneracy, i.e. $1+(1/p)<m<n/p$, see \cite{DA}. {\em It does not hold} in case $m=0$, see  Section \ref{contlp} below.   
Thus it remains to be proved in the following cases:
 \begin{itemize}
 \item $m=2$, $p=1$, and $n\geq 3$.
    \item $m=n/p>1$ and $p>1$.
     \item $m\geq \max(n,2)$ and $p=1$.
\end{itemize}
 If we except the space $W^2_1(\R)$, all the Sobolev spaces under consideration
  are particular cases of the Adams-Frazier spaces, or of the Sobolev algebras. We will prove the automatic continuity for those spaces,
and for their homogeneous counterparts, conveniently realized. Contrarily to the case $m=1$, where the proof of continuity of $T_f$ is much more difficult for $p=1$, see \cite[p.~219]{MM}, our proof in case $m\geq 2$ will cover all values of $p\geq 1$.

     \section*{Plan - Notation}
     
     In Section \ref{contlp} we recall the classical result on the continuity of $T_f$ in $L_p$ spaces. We take this opportunity to correct some
     erroneous statement in the literature. In Section \ref{AF}, we recall the characterization of composition operators acting in inhomogenous and homogeneous Adams-Frazier spaces, and in Sobolev algebras. In Section \ref{AFHom} we explain the specific difficulties concerning the continuity of $T_f$ in homogeneous spaces, which can be partially overcome by using realizations.
     Section \ref{Proof} is devoted to the proof of the continuity of $T_f$.\\
     
     We denote by ${\mathbb N}$ the set of all positive integers, including $0$. All functions occurring in the paper are assumed to be real valued. We denote by ${\mathcal P}_k$ the set of polynomials on $\re^n$, of degree less or equal to $k$. If $f$ is a function on $\re^n$, we denote by $[f]_k$ its equivalence class modulo ${\mathcal P}_k$. We consider a classical mollifiers sequence $\theta_\nu(x):= \nu^n \theta (\nu x)$, $\nu\geq 1$, where $\theta\in {\mathcal D}(\R)$ and $\int_{\R} \theta (x)\, \mathrm{d}x=1$. For all $N\in {\mathbb N}$, we denote by $C^N_b(\re)$ the space of functions $f:\re \rightarrow \re$, with continuous bounded derivatives up to order $N$.
 In all the paper, $m$ is an integer $>1$ and the real number $p$ satisfies $1\leq p <+\infty$. $W^m_p(\R)$ and  $\dot{W}^m_p(\R)$  are the classical inhomogeneous and homogeneous Sobolev spaces, endowed with the norms and seminorms
 \[\|g\|_{W^m_p}:= \sum_{|\alpha|\leq m} \|g^{(\alpha)}\|_p\,,\quad \|g\|_{\dot{W}^m_p}:= \sum_{|\alpha|= m} \|g^{(\alpha)}\|_p\,,\] respectively.
For topological spaces $E,F$, the symbol $E\hookrightarrow F$ means an imbedding, i.e.: $E\subseteq F$ and the natural mapping $E\rightarrow F$ is continuous.\\

The authors are grateful to Alano Ancona, Mihai Brancovan and Thierry Jeulin for fruitful discussions during the preparation of the paper.

   \section{ The case of $L_p$}\label{contlp}
   
In a recent survey paper on composition operator, the first author said that all the composition operators acting in $L_p(\R)$ are continuous
(see \cite{BS}, in particular the first line of the tabular at page $123$). {\em This assertion is erroneous}. Indeed $T_f$ takes $L_p(\R)$ to itself iff
 $|f(t)|\leq c\,|t|$ for some constant $c$, see  \cite[thm.~3.1]{AZ}. Clearly this property does not imply the continuity of $f$ outside of $0$.
 Instead we have the following:
 
  \begin{prop}\label{neccont} Let $(X,\mu)$ be a measure space. Let $f:\re\rightarrow \re$ be  s.t. $T_f$ takes $L_p(X,\mu)$ to itself. Then $T_f$ is continuous from $L_p(X,\mu)$ to itself iff $f$ is continuous.
\end{prop}

{\bf Proof.}  Let us assume that $T_f$ is continuous on $L_p$.
Without loss of generality, assume that $f(0)=0$. Let $A$ be  a measurable set in $X$ s.t. $0<\mu(A)<+\infty$. For all real numbers $u,v$, it holds
\be\label{fuv} \|f\circ  u\chi_A - f\circ v\chi_A\|_p= |f(u)-f(v)| \,\mu(A)^{1/p}\,.\ee
Clearly $v\to u$  in $\re$ implies $ v\chi_A \to u\chi_A$ in $L_p$.
By the continuity of $T_f$, and by (\ref{fuv}), we obtain the continuity of $f$. For the reverse implication,
we refer to \cite[thm.~3.7]{AZ}. We can also use the following statement:

 \begin{prop}\label{contlpq} Assume $q\in [1,+\infty[$. Let $(X,\mu)$ be a measure space. Let $f:\re\rightarrow \re$ be a continuous function s.t.
 for some constant
$c>0$, it holds $|f(t)|\leq c\,|t|^{p/q}$, for all $t\in \re$.  Then $T_f$ is a continuous mapping from $L_p(X,\mu)$ to $L_q(X,\mu)$.\end{prop}

{\bf Proof.} We follow \cite[thm.~2.2]{AP}. Let $(g_\nu)$ be a sequence converging to $g$ in $L_p(X,\mu)$.
By a classical measure theoretic result\footnote{See the proof of the completeness of $L_p$, e.g. \cite[thm.~3.11, p.~67]{R2}}, there exists a subsequence  $(g_{\nu_k})$ and a function $
h\in L_p(X,\mu)$ s.t. 
\[g_{\nu_k} \rightarrow g\quad \mathrm{a.e.} \,,\quad \left|g_{\nu_k} \right| \leq h\,.\]
By the continuity of $f$, it holds $f\circ g_{\nu_k} \rightarrow  f\circ g$ a.e..
By assumption on $f$, it holds
\[ \left |f\circ g_{\nu_k} -f\circ g\right|^q\leq (2c)^q\,h^p\,.\]
By Dominated Convergence Theorem, we conclude that $\left \|f\circ g_{\nu_k} -f\circ g\right\|_q$ tends to $0$.

  \section{Adams-Frazier spaces and related spaces}\label{AF}

\subsection{Function spaces}

The  inhomogeneous and homogeneous Adams-Frazier spaces are defined as follows:
\[A^m_p(\re^n):=W^{m}_{p}\cap  \dot{W}^{1}_{mp}(\re^{n})\,,\quad \dot{A}^m_p(\re^n):=\dot{W}^{m}_{p}\cap  \dot{W}^{1}_{mp}(\re^{n})\,.\]
Both spaces are endowed with their natural norm and seminorm:
\[ \|f\|_{A^m_p}:= \| f\|_{W^{m}_{p}}+ \|f\|_{ \dot{W}^{1}_{mp}}\,,\quad
 \|f\|_{\dot{A}^m_p}:= \|f\|_{\dot{W}^{m}_{p}}+ \| f\|_{\dot{ W}^{1}_{mp}}\,.\]
 The pertinency of those spaces w.r.t. composition operators was first noticed in \cite{AF}, see also the Introduction of \cite{Bou_09a}.
 By Sobolev imbedding, it holds  $W^m_p(\R)\hookrightarrow W^1_{mp}(\R)$, in case $m\geq n/p$, hence 
 \be\label{surcrit} m\geq n/p\quad \Rightarrow \quad A^m_p(\re^n)=W^{m}_{p}(\R)\,.\ee
 In particular the critical Sobolev spaces $W^m_p(\R)$, $m=n/p$, are Adams-Frazier spaces.\\
 
  The intersections $\dot{W}^m_p\cap L_\infty (\R)$  and $W^m_p\cap L_\infty (\R)$ are subalgebras of $L_\infty(\R)$ for the usual pointwise product, see Remark \ref{algebra} below. We call them the homogeneous  and inhomogeneous {\em Sobolev algebras}, and we endow them with their natural norms. Let us notice that the second is a proper subspace of the first, since the nonzero constant functions do not belong to $W^m_p(\R)$. By the Gagliardo-Nirenberg inequalities, see e.g. \cite[(6), p.~6108]{Bou_09a}, we have the imbeddings 
  \be\label{SAtoAF} \dot{W}^m_p\cap L_\infty (\R)\hookrightarrow \dot{A}^m_p(\R)\,,\quad W^m_p\cap L_\infty (\R)\hookrightarrow A^m_p(\R)\,.\ee
  In particular $W^m_p(\R)$ coincides with the corresponding Sobolev algebra if $m>n/p$, or $m=n$ and $p=1$.
   The following statement characterizes the Adams-Frazier spaces which coincide with the corresponding Sobolev algebras:
  
   \begin{prop}\label{embedlinfty} 
   1- The inclusion $A^m_p(\re^n)\subset L_\infty (\R)$ holds iff $m>n/p$, or $m=n$ and $p=1$.
   
   2- The inclusion $\dot{A}^m_p(\re^n)\subset L_\infty (\R)$ holds iff $m=n$ and $p=1$.
   \end{prop}
   
    {\bf Proof.} 1- In case $m\geq n/p$, it suffices to apply (\ref{surcrit}) and the known properties of Sobolev spaces.
    Now assume $m<n/p$. Define $f_\lambda(x):= |x|^\lambda \varphi (x)$, where $\varphi\in {\mathcal D}(\R)$ and $\varphi (x)=1$ near $0$.
    Then $f_\lambda \in W^m_p(\R)$ iff $\lambda > m-(n/p)$. If we take
    \[ 1- \frac{n}{mp}<\lambda <0\,,\]
    we obtain $f_\lambda \in A^m_p(\re^n)$ but $f_\lambda\notin L_\infty (\R)$.
    
   2- In case $m<n/p$, or $m=n/p$ and $p>1$, the first part implies a fortiori $\dot{A}^m_p(\re^n)\not\subset L_\infty (\R)$.
    
     Now assume $m>n/p$. Define $f_\lambda(x):= |x|^\lambda (1- \varphi (x))$, with the same function $\varphi$ as above.
    Then $f_\lambda \in \dot{W}^m_p(\R)$ iff $\lambda <m-(n/p)$. If we take
    \[ 0<\lambda <1- \frac{n}{mp}\,,\]
    we obtain $f_\lambda \in \dot{A}^m_p(\re^n)$ but $f_\lambda\notin L_\infty (\R)$.
    
   If $f\in \dot{W}^n_1(\R)$, there exists $g\in C_0(\R)$ s.t.
   $f-g\in {\mathcal P}_{n-1}$, see \cite[thm.~3]{111-1}. If moreover $f\in \dot{W}^1_n(\R)$, one proves easily that $f-g$ is a constant.
Thus we obtain the inclusion  $\dot{A}^n_1(\re^n)\subset C_b(\R)$. 

  \begin{rem}\label{w21notaf} {\em The above proof shows also that $W^2_1(\R)$ does not coincide with $A^2_1(\R)$ in case $n\geq 3$: if we consider the function $f_\lambda$ of the above first part, with
  \[ 2-n< \lambda < 1- \frac{n}{2}\,,\]
  then $f_\lambda\in W^2_1(\R)$, but $f_\lambda\notin \dot{W}^1_2(\R)$.}\end{rem}

  \subsection{Uniform localization}
Let us introduce a function $\psi\in \mathcal{D}(\re)$, positive,  s.t.
  \[ \sum_{\ell \in \Z} \psi(x-\ell) = 1\quad \mathrm{for\,all}\,\, x\in \re\,.\]
  If $f$ is any distribution in $\re$, we set $f_\ell(x):=f(x)\psi(x-\ell)$, and we call the decomposition
  $f= \sum_{\ell \in \Z}f_\ell$ a {\em locally uniform decomposition} (abbreviated  as LU-decomposition) of $f$.
  If $E$ is a normed space of distributions in $\re$, we denote by $E_{lu}$ the set of distributions $f$ such that
  \[\|f\|_{E_{lu}}:= \sup_{\ell\in \Z} \|f_\ell\|_{E} <+\infty\,.\]
  Under standard assumptions on $E$, it is known that the space $E_{lu}$ does not depend on the choice of the function $\psi$.
 We refer to \cite{Bou_88-2,116-1,117-1} for details. In case $E=L_p(\re)$, $E_{lu}$ is denoted by $L_{p,lu}(\re)$.  
   For $N\in \N$, let us denote by  $\dot{W}^{N}_{L_{p,lu}}(\re)$ the homogeneous Sobolev space based upon  $ L_{p,lu}(\re)$, i.e. the set of functions $f$ s.t. $f^{(N)}\in  L_{p,lu}(\re)$. The inhomogeneous Sobolev space ${W}^{N}_{L_{p,lu}}(\re)$ is defined similarly. 
  Both spaces are endowed with their natural norm and seminorm:
 \[  \|f\|_{ \dot{W}^{N}_{L_{p,lu}}}:= \|f^{(N)}\|_{ L_{p,lu}}\,,\quad  \|f\|_{W^{N}_{L_{p,lu}}}:= \sum_{j=0}^N\|f^{(j)}\|_{ L_{p,lu}}.\]

 Let us recall that ${W}^{N}_{L_{p,lu}}(\re)$ coincides with $(W^N_p(\re))_{lu}$, with equivalent norms. As a consequence, we have the following property:
 
  \begin{prop}\label{anyI} Let $N\in \N$. For all compact interval $I$ of $\re$, there exists a constant $c=c(N,p,I)>0$ s.t.
\[\left\| \sum_{\ell\in \Z} g_\ell\right\|_{W^N_{L_{p,lu}}} \leq c\, \sup_{\ell\in\Z} \|g_\ell\|_{W^N_p}\,,\]
for all sequence $(g_\ell)_{\ell\in \Z}$ s.t. the function $g_\ell$ is supported by the translated interval $I+\ell$, for all $\ell\in \Z$.
\end{prop}

{\bf Proof.} See \cite[lem.1]{Bou_88-2}.

  \begin{prop}\label{cinfty} For any $N\in \N$, $C^\infty\cap W^N_{L_{p,lu}}(\re)$ is a dense subspace of  $W^N_{L_{p,lu}}(\re)$.\end{prop}
  {\bf Proof.}   We use the mollifiers sequence $(\theta_\nu)_{\nu\geq 1}$ introduced in Notation, with $n=1$. Let $f\in W^N_{L_{p,lu}}(\re)$ and $\varepsilon>0$. Let $f=\sum_{\ell \in \Z}f_\ell$ be the LU-decomposition of $f$. Let us set
 \[v(x):= \sum_{\ell\in \Z} \theta_{j_\ell}\ast f_\ell\,,\]
for some convenient sequence
 $(j_\ell)_{\ell\in \Z}$. First we observe that the sum which defines the function $v$ is locally finite.
 Indeed $\theta_{j_\ell}\ast f_\ell$ is supported by $I+\ell$, for some fixed compact interval $I$. As a consequence,
 $v$ is a $C^\infty$ function. By condition $p<\infty$ and by classical properties of Sobolev spaces, we can choose $j_\ell$ s.t.
 $ \|\theta_{j_\ell}\ast f_\ell-f_\ell\|_{W^N_p}\leq \varepsilon$. By Proposition \ref{anyI}, we deduce
 \[\|v-f\|_{W^{N}_{L_{p,lu}}}\leq c\,\sup_{\ell\in \Z} \|\theta_{j_\ell}\ast f_\ell-f_\ell\|_{W^N_p}\leq c\,\varepsilon\,.\]
This ends up the proof.

   \begin{prop}\label{alt} For any integer $N>0$, it holds $W^{N}_{L_{p,lu}}(\re)=\dot{W}^{N}_{L_{p,lu}}\cap L _{\infty}(\re)\hookrightarrow C^{N-1}_{b}(\re)$.\end{prop}

{\bf Proof.}  1- Let $f= \sum_{\ell\in \Z} f_\ell$ be the LU-decomposition of $f \in W^{N}_{L_{p,lu}}(\re)$. By condition $N\geq 1$, we have the Sobolev imbedding $W^N_p(\re)\hookrightarrow C_b(\re)$, hence
\[ \sup_{\ell\in \Z} \|f_\ell\|_{\infty} \leq c_1\,\sup_{\ell\in \Z} \|f_\ell\|_{W^N_p} \leq c_2 \,\|f\|_{W^{N}_{L_{p,lu}}}\,.\]
By considering the support of $f_\ell$, we conclude that $f\in C_b(\re)$.

2- Let $f\in\dot{W}^{N}_{L_{p,lu}}\cap L _{\infty}(\re)$. By \cite[lem.~1]{Bou_91}, it holds $f \in W^{N-1}_{\infty}(\re)$.
A fortiori, it holds $f^{(j)}\in L_{p,lu}(\re)$ for all $j=0,\ldots, N-1$. Hence $f\in W^{N}_{L_{p,lu}}(\re)$.

3-  If $f\in W^{N}_{L_{p,lu}}(\re)$, then $f\in W^{N}_{p}(\re)_{loc}$. By Sobolev imbedding, it follows that $f\in C^{N-1}_b(\re)$.
\subsection{Composition operators in Adams-Frazier spaces}

 \begin{prop}\label{compinhom}  The composition operator $T_f$ takes ${A}^m_p(\re^n)$ to itself if
 $f'\in  W^{m-1}_{L_{p,lu}}(\re)$ and $f(0)=0$. Moreover, for some constant $c=c(n,m,p)>0$, it holds
  \be\label{esticompo}
 \|f\circ g\|_{A^{m}_p} \leq c\,  \|f'\|_{ W^{m-1}_{L_{p,lu}}} \left(1+ \|g\|_{A^{m}_{p}} \right)^m \,,\ee
 for all $g\in A^m_p(\re^n)$. 
   \end{prop}

 \begin{prop}\label{comphomo}  The composition operator $T_f$ takes $\dot{A}^m_p(\re^n)$ to itself if
 $f'\in  W^{m-1}_{L_{p,lu}}(\re)$.   Moreover, the estimation (\ref{esticompo}) holds true with $A$ replaced by $\dot{A}$.
 \end{prop}
 
 We refer to \cite[thm.~1, p.~6107]{Bou_09a}, \cite[thm.~25]{BS} and to Proposition \ref{alt}. 
 The condition  $f'\in  W^{m-1}_{L_{p,lu}}(\re)$ is not only sufficient, but also necessary, for all Adams-Frazier spaces which are {\em not included into}
 $L_\infty(\R)$. 
 For such spaces, Propositions \ref{compinhom} and \ref{comphomo}
  constitute characterizations of the  acting composition operators.
  For spaces imbedded into $L_\infty(\R)$, we have the following alternative result:  
 
  \begin{thm}\label{compalg} 
  The operator $T_f$ takes $\dot{W}^m_p\cap L_\infty (\R)$ to itself iff
 $f\in  W^{m}_p(\re)_{loc}$, and 
    $W^m_p\cap L_\infty (\R)$ to itself iff
 $f\in  W^{m}_p(\re)_{loc}$ and $f(0)=0$.
 \end{thm}
 
 {\bf Proof.}   This statement is classical, see \cite{Bou_91} and \cite[thm.~2]{BS}. We recall here the part of the proof which will be useful to prove the continuity of $T_f$.
  For every $r>0$, we denote by $B_r$ the ball of center $0$ and radius $r$ in $\dot{W}^m_p\cap L_\infty (\R)$.
  It will suffice to prove that {\em $T_f$ takes $B_r$ to $\dot{W}^m_p\cap L_\infty (\R)$,} for all $r>0$.
We introduce a family of auxiliary functions $\omega_r\in \mathcal{D}(\re)$ s.t. $\omega_r (t)=1$ for $|t|\leq r$. 
Then
  \be\label{inball} \forall g\in B_r\,,\qquad T_f(g)= T_{f\omega_r}(g)\,.\ee
  If $f\in  W^{m}_p(\re)_{loc}$, it is clear that $(f\omega_r)'\in W^{m-1}_{L_{p,lu}}(\re)$.
We can apply Proposition \ref{comphomo} to $f\omega_r$.  By (\ref{SAtoAF}), we conclude that $T_{f\omega_r}$ takes $\dot{W}^m_p\cap L_\infty (\R)$ to $\dot{W}^m_p(\R)$. By condition $m\geq 1$, the Sobolev imbedding yields $f\omega_r\in C_b(\re)$, and we obtain that $T_{f\omega_r}$ takes $\dot{W}^m_p\cap L_\infty (\R)$ to $L_\infty (\R)$. This ends up the proof.

 \begin{rem}\label{algebra} {\em
 In particular, any function of class $C^m$ acts on $\dot{W}^m_p\cap L_\infty (\R)$ by composition. Applying this to the function $f(t):=t^2$, we derive
 immediately the algebra property.} \end{rem}

For $n=1,2$, the space $W^2_1(\R)$ is a Sobolev algebra, for which the acting composition operators are
described in Theorem \ref{compalg}. In the other cases, we have the following result (see \cite{Bou_91}):

\begin{thm}\label{compw21} Assume that 
 $n\geq 3$. Then $T_f$ takes $W^2_1(\R)$ to itself iff $f(0)=0$ and $f''\in L_1(\re)$. Moreover there exists $c=c(n)>0$ s.t.
     \be \label{estw21}  \|f\circ g\|_{W^2_1} \leq c\,(|f'(0)|+ \|f''\|_1)\,   \|g\|_{W^2_1}\,,\ee
     for all such $f$'s and all $g\in W^2_1(\R)$.
     \end{thm}

  \section{ Homogeneous spaces and their realizations}\label{AFHom}

Usually, an homogeneous function space $F$, such as $\dot{W}^m_p(\R)$, is only a seminormed space, with $\|f\|=0$ iff $f\in{\mathcal P}_k$, for some $k\in \mathbb{N}$ depending on $F$. The presence of polynomials, with a seminorm equal to $0$, has some pathological effects on composition operators. Recall, for instance, 
the following (see \cite[prop.~11]{Bou_09a}):

\begin{prop}\label{triv} If $m>1$ and $n>1$,
the only functions $f$, for which $T_f$ takes $\dot{W}^m_p(\R)$ to itself, are the affine ones.\end{prop}

This degeneracy phenomenon does not occur in homogeneous Adams-Frazier spaces, see Proposition \ref{comphomo}. However, the continuity of
$T_f$ is a tricky question. The statement: `` $T_f$ is continuous as a mapping of the seminormed space $\dot{A}^m_p(\re^n)$ to itself '' makes sense, but it has no chance to be true. Assume that, for a sequence $(g_\nu)$ tending to $g$ in  $\dot{A}^m_p(\re^n)$,  the sequence $(f\circ g_\nu)$ tends to $f\circ g$ in  $\dot{A}^m_p(\re^n)$. Then, for a sequence $(c_\nu)$ of real numbers, the sequence $(g_\nu+c_\nu)$ tends also to $g$ in $\dot{A}^m_p(\re^n)$. But the sequence
$(f\circ (g_\nu+c_\nu)- f\circ g)$ cannot tend to $0$ in $\dot{A}^m_p(\re^n)$, {\em whatever be} the sequence $(c_\nu)$. If, for instance, $f(t):=\sin t$,  $g$ is a nonzero function in ${\mathcal D}(\R)$, and $g_\nu:=g$ for all $\nu\in {\mathbb N}$, then
$f\circ (g+\pi)- f\circ g=-2 f\circ g$, a function which is not constant. \\

In order to avoid the disturbing effect of polynomials, two ideas seems available.
The first one would be to consider the factor space  $F/{\mathcal P}_k$. But that does not work.  Indeed, if $g_1$ and $g_2$ differ by a polynomial, the same does not hold for $f\circ g_1$ and $f\circ g_2$, hence we cannot extend the operator $T_f$ to the factor space. 
The second one consists in restricting $T_f$ to a
vector subspace $E$ s.t. $
F= E \oplus \mathcal{P}_k$.
We will exploit this idea in case of Adams-Frazier spaces  $\dot{A}^m_p(\re^n)$, and the space $\dot{W}^2_1(\R)$.

\subsection{Realizations of homogeneous Adams-Frazier spaces}

Let us begin with a remark: 
\begin{prop}\label{banach} \begin{enumerate}
\item The factor space $\dot{A}^m_p(\R)/\mathcal{P}_0$, endowed with the {\em norm} $\|-\|_{\dot{A}^m_p}$, is a Banach space. 
\item Any subspace $E$ of $\dot{A}^m_p(\R)$ s.t.
\be\label{realAF} \dot{A}^m_p(\R)= E \oplus \mathcal{P}_0\ee
 is a Banach space for the norm $\|-\|_{\dot{A}^m_p}$.
 \end{enumerate}
\end{prop}

{\bf Proof.} 1- It is well known that $\dot{W}^m_p(\re^n)/\mathcal{P}_{m-1}$ is a Banach space if endowed with the norm $\|-\|_{\dot{W}^m_p}$, see \cite[1.1.13, thm.~1]{MA}.
If $([g_\nu]_0)$ is a Cauchy sequence in $\dot{A}^m_p(\re^n)/\mathcal{P}_{0}$, then there exist $u\in \dot{W}^m_p(\re^n)$ and $v\in \dot{W}^1_{mp}(\re^n)$ s.t. $(g_\nu)$ tends to $u$ in $ \dot{W}^m_p(\re^n)$ and to $v$ in $ \dot{W}^1_{mp}(\re^n)$. By \cite[1.1.13, thm.~2]{MA}, there exist a sequence $(r_\nu)$ in $\mathcal{P}_{m-1}$ and a sequence $(c_\nu)$ in $\mathcal{P}_{0}$, s.t. $(g_\nu -r_\nu)$ tends to $u$, and $(g_\nu-c_\nu)$ to $v$, in $ L_{1}(\re^n)_{loc}$.
Thus $(r_\nu -c_\nu)$ is a sequence in $\mathcal{P}_{m-1}$, which tends to $v-u$ in $ L_{1}(\re^n)_{loc}$. We conclude that $v-u\in \mathcal{P}_{m-1}$, and that $([g_\nu]_0)$ tends to $[v]_0$ in $\dot{A}^m_p(\re^n)/\mathcal{P}_0$. 

2- If $E$ is a subspace of $\dot{A}^m_p(\re^n)$ s.t. (\ref{realAF}), the linear map $f\mapsto [f]_0$ is an isometry from $E$  onto $\dot{A}^m_p(\re^n)/\mathcal{P}_{0}$. The completeness of $E$ follows by the first part. This ends up the proof.\\

A subspace $E$ satisfying (\ref{realAF}) will be of interest only if it is a Banach space of distributions. This motivates the following definition:

\begin{DE}\label{honest} A subspace $E$ of $\dot{A}^m_p(\R)$ s.t.
(\ref{realAF}) is called
a {\em realization} of $\dot{A}^m_p(\re^n)$ if one of the following equivalent properties holds:
\begin{enumerate}
\item the inclusion mapping $E\rightarrow \mathcal{S}'(\R)$ is continuous;
\item the inclusion mapping $E\rightarrow  L_{1}(\R)_{loc}$ is continuous;
\item for every sequence $(g_\nu)$ tending to $g$ in  $E$, there exists a subsequence $(g_{\nu_k})$ s.t. $g_{\nu_k}\rightarrow g$ a.e..
\end{enumerate}
\end{DE}

The equivalence between the three properties follows easily by the Closed Graph Theorem.

\begin{rem}\label{realequiv} {\em In \cite{Bou_09a}, we used a slightly weaker definition for a realization of $\dot{W}^m_p(\R)$.
We said that a subspace $E$ of $\dot{W}^m_p(\R)$ is a realization if
\be\label{oldreal1} \dot{W}^m_p(\R) = E\oplus \mathcal{P}_{m-1}\,.\ee
If (\ref{oldreal1}) holds we obtain a linear mapping
$\sigma: \dot{W}^m_p(\R)/\mathcal{P}_{m-1} \rightarrow \mathcal{S}'(\R)$ s.t.
\[ \forall u\in \dot{W}^m_p(\R)/\mathcal{P}_{m-1}\quad  [\sigma (u)]_{m-1} = u\,,\]
and whose range is $E$.
Then $\sigma$ is a realization, in the sense of  \cite{088-1,111-1,113-1}, if $\sigma$ is a continuous
mapping from $\dot{W}^m_p(\R)/\mathcal{P}_{m-1}$ to
 $\mathcal{S}'(\R)$: this is precisely what means Definition \ref{honest}.

 }\end{rem}

Now we turn to the description of the usual realizations of $\dot{A}^m_p(\R)$.
Except in case $m=n$, $p=1$, it will suffice to realize $\dot{W}^1_{mp}(\re^n)$, then restrict to $\dot{A}^m_p(\re^n)$.
The most natural realizations are those which retain the invariance properties of  $\dot{A}^m_p(\re^n)$ w.r.t. translations or dilations.
It is classically known that such realizations do not always exist, see \cite{088-1,111-1,113-1}.\\

{\em 1- Case $m<n/p$.}  Let us set
\[\frac1q:= \frac{1}{mp}-\frac1n\,.\]
Then $L_q\cap \dot{W}^1_{mp}(\re^n)$ is a realization of $\dot{W}^1_{mp}(\re^n)$, see \cite[prop.~14]{Bou_09a}. Hence
$L_q\cap \dot{A}^m_p(\re^n)$ is a realization of $ \dot{A}^m_p(\re^n)$. Clearly it is invariant w.r.t. translations and dilations.\\

 {\em 2- Case $m>n/p$.} By condition $1> \frac{n}{mp}$,  $\dot{W}^1_{mp}(\re^n)$ is a subset of $C(\R)$. Then the subspace
 $\{ f\in \dot{A}^m_p(\re^n)\,:\, f(0)=0\}$ is a
 dilation invariant realization of $\dot{A}^m_p(\re^n)$.\\
 
 {\em 3- Case $m=n$ and $p=1$.} As observed in the proof of Proposition \ref{embedlinfty}, $C_0\cap \dot{A}^n_1(\re^n)$
 is a realization of $\dot{A}^n_1(\re^n)$, clearly invariant w.r.t. translations and dilations.\\

 {\em 4- Case $m=n/p$ and $p>1$.} In such a case, $\dot{A}^m_p(\R)$ does not admit invariant realizations. This  can be deduced from \cite[thms.~5.4, 5.7]{113-1}.\\
 
 In all cases we can use ``rough'' realizations described as follows. Let us recall that $\dot{A}^m_p(\re^n)$ is a subset of $L_{pm}(\R)_{loc}$, see \cite[1.1.2]{MA}. Let us define $q$ by
\[\frac1q := 1- \frac{1}{mp}\,.\]
Let $\varphi$ be a compactly supported measurable function  in $L_q(\re^n)$, s.t.
\[ \int_{\re^n} \varphi (x)\,{\rm d}x =1\,.\]
We can define a linear functional on  $\dot{A}^m_p(\re^n)$ by setting 
\[\Lambda (g):= \int_{\R} \varphi (x)g(x)\,{\rm d}x\,.\]
 Then the kernel of $\Lambda$ is a realization of $\dot{A}^m_p(\re^n)$, with no invariance property.

\subsection{Realizations of $\dot{W}^2_1(\R)$}\label{realw21}

According to Proposition \ref{triv}, there is no nontrivial composition operator which takes $\dot{W}^2_1(\R)$ to itself if $n>1$. In such a case, we are forced to introduce realizations,
i.e. subspaces  $E$ s.t.
$ \dot{W}^2_1(\R)= E \oplus \mathcal{P}_1$, and satisfying the equivalent properties of Definition \ref{honest}.
Let us recall the known results concerning invariant realizations, and
composition operators acting on them, see \cite{111-1} and \cite[prop.~18]{Bou_09a} for details.\\

{\em 1- Case $n=1$.}  $T_f$ takes $\dot{W}^2_1(\re)$ to itself iff $f\in \dot{W}^2_1(\re)$, and, for a such $f$, it holds
\be\label{n=1} 
\|(f\circ g)''\|_1 \leq c\, \left( \|f''\|_1 + |f'(0)|\right)\, \left( \|g''\|_1 + |g'(0)|\right)\,,
\ee
for every $g\in \dot{W}^2_1(\re)$. 

\begin{LEM}\label{honestn=1}  For all ${\bf{\alpha}}:=(\alpha_1,\alpha_2,\alpha_3)\in \re^3$ s.t. $\alpha_1+\alpha_2+\alpha_3=1$, the subspace
\be\label{realn=1} E_{ \alpha}= \{ g\in \dot{W}^2_1(\re)\,:\, g(0)=0\,,\,\, \alpha_1 g'(-\infty)+ \alpha_2 g'(0)+ \alpha_3 g'(+\infty)=0\}\ee
is a dilation invariant realization of $\dot{W}^2_1(\re)$. 
\end{LEM}

{\bf Proof.} Let us recall that {\em every function in $\dot{W}^1_1(\re)$ is continuous, with finite limits at $\pm\infty$.}
Thus the definition of $E_\alpha$ makes sense. It is easily seen that $ \dot{W}^2_1(\re)= E_\alpha \oplus \mathcal{P}_1$.
If $g\in \dot{W}^2_1(\re)$ and $g(0)=0$, it holds
$\|g'\|_\infty\leq | g'(-\infty)|+ \|g''\|_1$, hence $|g(x)|\leq |x|\,(| g'(-\infty)|+  \|g''\|_1)$ for all $x\in \re$.
It follows that $E_\alpha$ is imbedded into $L_{1}(\R)_{loc}$.


\begin{rem}{\em It can be proved that
 any dilation invariant realization of $\dot{W}^2_1(\re)$ is necessarily equal to  $E_{\bf \alpha}$ for some ${\bf \alpha}$, see  \cite[prop.~11]{111-1}.}
 \end{rem}

{\em 2- Case $n=2$.} According to \cite[thm.~3]{111-1}, 
\be\label{realn=2}E:=C_0\cap \dot{W}^2_1(\re^2)\ee is a realization of  $\dot{W}^2_1(\re^2)$. Indeed, it is the unique translation invariant realization, see \cite[thm.~6]{111-1}. By Theorem \ref{compalg}, $T_f$ takes $E$ to $\dot{W}^2_1(\re^2)$ iff $f\in W^2_1(\re)_{loc}$.\\

{\em 3- Case $n\geq 3$.} According to  \cite[thm.~2]{111-1}, if $\frac1q := 1- \frac2n$, then
\be\label{realn>2} E:=L_q\cap \dot{W}^2_1(\R)\ee is a realization of  $\dot{W}^2_1(\R)$. Indeed, it is the unique translation invariant realization, and the unique dilation invariant realization, of $\dot{W}^2_1(\R)$, see
 \cite[thm.~6, prop.~11]{111-1}. $T_f$ takes $E$ to $\dot{W}^2_1(\re)$ iff $f\in \dot{W}^2_1(\re)$. The estimation
  \be \label{esthomw21}  \|f\circ g\|_{\dot{W}^2_1} \leq c\, (|f'(0)|+ \|f\|_{\dot{W}^2_1})\,  \|g\|_{\dot{W}^2_1}\ee
   holds  for all $f\in \dot{W}^2_1(\re)$ and all $g\in E$.

\section{Continuity theorems}\label{Proof}


\begin{thm}\label{contAF} Let $f:\re\rightarrow \re$ be s.t. $f' \in W^{m-1}_{L_{p,lu}}(\re)$. Let $E$ be a realization of $\dot{A}^m_p(\R)$. 
Then $T_f$ is continuous from $E$ to $\dot{A}^m_p(\R)$.
If moreover $f(0)=0$, then $T_f$
is continuous from $A^m_p(\R)$  to itself.
\end{thm}

Under the stronger assumption $f' \in C^{m-1}_b(\re)$, the continuity of $T_f$ on $A^m_p(\R)$ is a classical result,
seemingly with no reference in the literature;  in their article on composition operators in fractional Sobolev spaces \cite{BM},
 Brezis and Mironescu  said only that the proof is ``very easy via the standard Gagliardo-Nirenberg inequality''.
 This ``folkloric'' proof will be recalled below, see the proof of Lemma \ref{low}:  what we do for terms with $s<m$ works as well for $s=m$, in case $f' \in C^{m-1}_b(\re)$.

   \begin{thm}\label{contalg} Let $f\in  W^{m}_p(\re)_{loc}$. Then $T_f$ is continuous from $\dot{W}^m_p\cap L_\infty (\R)$ 
 to itself. If moreover $f(0)=0$, then $T_f$ is continuous from $W^m_p\cap L_\infty (\R)$
  to itself.
  \end{thm}
  
  In case $p>1$, Theorem \ref{contalg} has been proved in \cite[cor.~2]{BL}, as a particular case of a continuity theorem for composition
  in Lizorkin-Triebel spaces. F.Isaia has also proved it for $W^m_p(\R)$, in case $m>n/p$ and $p\geq 1$, but with a stronger condition on $f$, namely 
  $f\in  W^{m}_\infty(\re)_{loc}$, see \cite[thm.~2.1, (iii)]{IS}.
  
   \begin{thm}\label{contw21}  Let $E$ be the realization of  $\dot{W}^2_1(\R)$ defined by
   (\ref{realn=1}) or   (\ref{realn=2}) or (\ref{realn>2}) according to the value of $n$.
 Let $f:\re\rightarrow \re$ be s.t. $T_f$ takes $E$ to  $\dot{W}^2_1(\R)$. Then
 $T_f$ is continuous from $E$ to  $\dot{W}^2_1(\R)$. If moreover $f(0)=0$, then $T_f$
is continuous from $W^2_1(\R)$  to itself.
     \end{thm}
     
      Let us notice that Theorem \ref{contw21} is less general than
     Theorem \ref{contAF} since we do not consider all the  realizations, but only the invariant ones.

 \subsection{Main tools}
 
Four propositions will be useful, where the first is elementary, the second follows by H\"older inequality and the third is classical, see e.g. \cite[I \S8.2, prop.~2 , p.~100]{baki_1994}.

\begin{prop}\label{t^alpha} If $\alpha>1$, the functions
 $t\mapsto |t|^\alpha$ and $t\mapsto \mathrm{sgn}(t)\,|t|^\alpha$ are  of class $C^1$ on $\re$.
 \end{prop}

  \begin{prop}\label{multili} Let  $p_1,\ldots , p_s\in [1,+\infty[$ s.t.
  \[ \sum_{j=1}^s \frac{1}{p_j} = \frac{1}{p}\,.\]
  Then the mapping $(\varphi_1,\ldots, \varphi_s)\mapsto \varphi_1\, \cdots
 \varphi_s$ is continuous from  $L_{p_1}\times\cdots \times L_{p_s}$ to $L_p$.\end{prop}

  \begin{prop}\label{algelem} Let us denote $I:= \{1,\dots , m\}$. Then, for every $(x_1,\ldots ,x_m)\in \re^m$, it holds
 \[ x_1x_2\cdots x_m=\frac{(-1)^m}{m!}\sum_{H\subseteq I}(-1)^{|H| }\left( \sum_{k\in H}x_k\right)^m \,,\]
 where the sum runs over all the subsets of $I$, and $|H|$ denotes the cardinal of $H$.
 \end{prop}

 \begin{prop}\label{gencont} Let $E$ be a Banach space of distributions in $\R$ s.t. $E\hookrightarrow L_{1}(\R)_{loc}$.
Let $T$ be a continuous mapping from $E$ to $L_p(\R)$. Let $\Phi \in C_b(\re)$. Define the mapping $V:E\rightarrow L_p(\R)$ by
$V(g):= (\Phi\circ g)\, T(g)$. Then $V$ is continuous from $E$ to $L_p(\R)$.
\end{prop}

{\bf Proof.} Let $(g_\nu)$ be a sequence s.t. $\lim g_\nu = g$ in $E$. By the imbedding $E\hookrightarrow L_{1}(\R)_{loc}$, we can assume that 
 $\lim g_\nu = g$ a.e., up to replacement by some subsequence.
 It holds
  \[ \|V(g_\nu) - V(g)\|_p\leq \|\Phi\|_\infty\, \|T(g_\nu) - T(g)\|_p + \left(\int_{\R} |\Phi\circ g_\nu -\Phi\circ g|^p\, |T(g)|^p\, {\rm d}x\right)^{1/p}\,.\]
 The second term of the above r.h.s. tends to $0$ by Dominated Convergence Theorem (Convergence a.e. follows by continuity of $\Phi$, domination by boundedness of $\Phi$ and the fact that $T(g)\in L_p(\R)$).\\

 \subsection{Proof of Theorem \ref{contAF}}

 In all the proof, we denote by $E$ the space $A^m_p(\R)$, or a realization of $\dot{A}^m_p(\R)$. 
Without loss of generality, we assume that $f$ is a {\em smooth} function s.t.
 $f' \in W^{m-1}_{L_{p,lu}}(\re)$, see Proposition \ref{cinfty}, and the estimation (\ref{esticompo}).
  
 \subsubsection{Continuity of $T_f$ from $E$ to $\dot{W}^1_{mp}$ and to $L_p$}
 
 By assumption $m\geq 2$, and by Proposition \ref{alt}, we have $f'\in C_b(\re)$. 
  For $j=1,\ldots ,n$, it holds
 $\partial_j(f\circ g)= (f'\circ g)\,\partial_jg$. 
Thus we can apply Proposition \ref{gencont} and conclude that
 $g\mapsto\partial_j(f\circ g)$ is a continuous mapping from $E$ to $L_{mp}(\R)$.
 We have also
 $\|f\circ g_1-f\circ g_2\|_p\leq \|f'\|_\infty \| g_1-g_2\|_{p}$, hence the continuity from $A^m_p(\R)$ to $L_{p}(\R)$, in case $f(0)=0$.
 
 \subsubsection{Continuity to $\dot{W}^m_{p}$: heart of the proof}\label{PREP}

Let us consider the nonlinear operator
\be\label{S_D} S_{D}(g) :=(f^{(m)}\circ g)\, (Dg)^m\,,\ee
where $D$ is any first order differential operator with constant coefficients, say
$D:= \sum_{j=1}^nc_j\partial_j$, for some real numbers $c_j$.  
The fact that $S_D$ takes $\dot{A}^m_p(\R)$ to $L_p(\R)$ is the heart of the proof of Proposition \ref{comphomo}, see \cite[par.~2.3]{Bou_09a}. Thus, we can expect that the following result will be the main argument for proving Theorem \ref{contAF}: 

\begin{LEM}\label{formula} The nonlinear operator $S_D$ is continuous from $E$ to $L_p(\R)$. \end{LEM}

{\bf Proof.}
Let  $u\in {\mathcal D}({\mathbb R})$ be such that
 $u\geq 0$ and
\[\forall y\in {\mathbb R}\,,\quad \sum_{\ell\in {\mathbb Z}} u^2(y-\ell) = 1\,.\]
Let us define
$\Phi_{\ell}(y):=u(y-\ell)\, | f^{(m)}
 (y)|^p$, for all $y\in \re$.
Since $f^{(m)}\in L_{p,lu}\cap C^\infty(\re)$, it holds 
 \[ \Phi_\ell \in C(\re) \quad \mathrm{and}\quad \sup_{\ell\in\mathbb{Z}}\|\Phi_{\ell} \|_1 < +\infty\,.\]
We define
  $\displaystyle{\Psi_{\ell}(y):= \int_{y}^{+\infty} \Phi_{\ell}(t)\,{\rm d}t}$.
Then 
  \be\label{linfty} \Psi_\ell \in C^1(\re) \quad \mathrm{and}\quad \sup_{\ell\in\mathbb{Z}}\|\Psi_{\ell} \|_\infty < +\infty\,.\ee
  
  Now we assume that $mp>2$ (the exceptional case, $m=2$ and $p=1$, needs a minor change, see below).
We can use Proposition \ref{t^alpha} with $\alpha:=mp-1$.
It holds
  \[ \|S_{D}(g)\|_p^p=\sum_{\ell \in {\mathbb Z}} \int_{\R} |f^{(m)}\circ g| ^p\, |Dg|^{mp}\, (u^2\circ (g-\ell))
 \,{\rm d}x\]
 \[ =\sum_{\ell \in {\mathbb Z}} \int_{\R} (\Phi_\ell\circ g)\, (Dg) \left( \mathrm{sgn}(Dg)\, |Dg|^{mp-1}\, (u\circ (g-\ell))\right)
\,{\rm d}x\,,\]
hence
\[ \|S_{D}(g)\|_p^p = \sum_{j=1}^n c_j\,\sum_{\ell \in{\mathbb Z}} V_{j,\ell}\,,\]
with
\[ V_{j,\ell}:= \int_{\R} (\Phi_\ell\circ g)\, \partial_jg\, \left( \mathrm{sgn}(Dg)\, |Dg|^{mp-1}\, (u\circ (g-\ell))\right)
\,{\rm d}x\,.\]
The computation of $V_{j,\ell}$ relies upon an integration by parts w.r.t the $j$-th coordinate.
This I.P. is justified by a classical theorem. {\em Any function which belongs locally to $W^1_p(\R)$
can be modified on a negligible set of $\R$, in such a way that the resulting function has the following property:
its restriction, on almost every line parallel to a coordinate axis, is locally absolutely continuous}. This theorem originates to Calkin \cite{CA}; see \cite[1.1.3, thm.1]{MA}
for the precise statement. Here, this theorem can be applied to the function $\Psi_\ell\circ g$, which belongs locally to $W^1_{mp}(\R)$, and to the function
\[w:= \mathrm{sgn}(Dg)\, |Dg|^{mp-1}\, (u\circ (g-\ell))\,,\]
which satisfies

\be\label{winw11} w\in \dot{W}^1_1\cap L_q(\R)\,,\quad q:= \frac{mp}{mp-1}\,.\ee

{\em Proof of (\ref{winw11}).}  Let $k=1,\ldots ,n$. By the assumption $g\in E$, and by the Gagliardo-Nirenberg inequality
 \[\|g\|_{\dot{W}^2_{mp/2}} \le c \|g\|_{\dot{W}^1_{mp}}^{\frac{m-2}{m-1}}\,\|g\|_{\dot{W}^m_{p}}^{\frac{1}{m-1}}\]
 (see \cite[(7), p.~6108]{Bou_09a}),
  it holds
$\partial_kDg\in L_{mp/2}(\R)$. On the other hand, since
\[|\partial_kw|\le c_1 |\partial_kDg|\, |Dg|^{mp-2}+c_2|Dg|^{mp-1}\,|\partial_kg|, \]
 the H\"older inequality applied to the r.h.s. (with exponents $\frac{mp}{2}$ and $\frac{mp}{mp-2}$ in the first term,   $\frac{mp}{mp-1}$ and $mp$ in the second)
gives $\partial_kw\in L_1(\R)$. \\

The property (\ref{winw11}) implies the following:\\

{\em For almost every $(x_1,\dots,x_{j-1}, x_{j+1},\ldots ,x_n)\in \re^{n-1}$, the function
\[t\mapsto w(x_1,\dots,x_{j-1}, t,x_{j+1},\ldots ,x_n)\] is absolutely continuous on $\re$, with limit $0$ at $\pm\infty$.}\\

Thus we obtain
 \[V_{j,\ell} = \int_{\R} (\Psi_\ell\circ g)\,\partial_j\left\{ \mathrm{sgn}(Dg)\, |Dg|^{mp-1}\, \, (u\circ (g-\ell))\right\}
\,{\rm d}x\,,\]
hence
 
\be\label{Sg} \|S_{D}(g)\|_p^p= \int_{\R} (F\circ g)  \, T_1(g)\,{\rm d}x
+  \int_{\R} (G\circ g) \,T_2(g)\, {\rm d}x \,,\ee
where
 \[F(y):= \sum_{\ell\in {\mathbb Z}} \Psi_\ell(y)\, u(y-\ell)\,,\quad
G(y):= \sum_{\ell\in {\mathbb Z}} \Psi_\ell(y)\, u'(y-\ell)\,,\]

and
  \[T_1(g):=  (mp-1) \,(D^2g)\, |Dg|^{mp-2}\,,\quad
T_2(g):=  |Dg|^{mp} \,.\]

In case $m=2$ and $p=1$, a similar computation starting with
 $$ \|S_{D}(g)\|_1= \int_{\R} |f''\circ g| \, (Dg)^{2}\,  \,{\rm d}x\,,$$ shows that (\ref{Sg}) is also valid.\\

 
    By Gagliardo-Nirenberg and by the definition of $E$, the linear mappings
   $g\mapsto D^2g$, $g\mapsto Dg$ are continuous from $E$ to $L_{mp/2}(\R)$ and $L_{mp}(\R)$ respectively. Thus, by using Proposition \ref{contlpq}, the mappings
    \[g\mapsto D^2g\,,\quad g\mapsto   |Dg|^{mp-2}\,,\quad g\mapsto   |Dg|^{mp}\] 
    are continuous from
  $E$ to $L_{\frac{mp}{2}}(\R)$, $L_{\frac{mp}{mp-2}}(\R)$ and $L_{1}(\R)$ respectively.
Since
\[ \frac{2}{mp} + \frac{mp-2}{mp} =1\,,\]
we can use Proposition \ref{multili}, and conclude that $T_1$ and $T_2$ are nonlinear continuous mappings from $E$ to $L_1(\R)$.  By the continuity of $\Psi_\ell$, and by (\ref{linfty}), it follows that $F$ and $G$ belong to $C_b(\re)$.
From (\ref{Sg}) and
Proposition \ref{gencont}, we deduce the following property:
  \be\label{contfunct} g\mapsto \|S_D(g)\|_p\,\, \mathrm{is\,\,a\,\,continuous\,\,nonlinear\,\,functional\,\,on}\,\, E\,.\ee
  Now consider a sequence $(g_\nu)$ which converges to $g$ in $E$. By the continuity of $f^{(m)}$, it holds
 \[\lim S_D(g_\nu) = S_D(g)\quad \mathrm{ a.e.}\,,\]
 up to replacement by a subsequence. By (\ref{contfunct}), it holds
  \[\lim \|S_D(g_\nu)\|_p = \|S_D(g)\|_p\,.\]
Then we can apply the Theorem of Scheff\'e \cite{SC}\footnote{Usually attributed to Scheff\'e, but first proved by F.~Riesz \cite{RI}, see the survey of N.~Kusolitsch \cite{KU}.} and conclude that
$\lim S_D(g_\nu) = S_D(g)$ in $L_p(\R)$.

 \subsubsection{Continuity to $\dot{W}^m_{p}$: end of the proof}\label{PREPII}
 
   
  We have to prove the continuity of $g\mapsto (f\circ g)^{(\alpha)}$, from $E
  $ to $L_p(\R)$, for all $\alpha \in {\mathbb N}^n$ s.t. $|\alpha|=m$.
We use the classical  formula
\begin{equation}\label{faa}
(f\circ g)^{(\alpha)}=\sum c_{\alpha,s, \gamma} (f^{(s)}\circ g)\, g^{(\gamma_{1})} \cdots
 g^{(\gamma_{s})}\,,\end{equation}
 where the parameters satisfy the conditions
 \begin{equation}\label{faa1}s=1,\ldots , m\,,\, |\gamma_{r}|>0\, (r=1,\dots ,s)\,,\, \sum_{r=1}^s\gamma_{r}= \alpha\,,    \end{equation}
 and the $c_{\alpha,s,\gamma}$'s are some combinatorial constants.
 Formula (\ref{faa}) holds true for any smooth function $g$. We need to prove it for any function $g\in E$.\footnote{In \cite[par.~2.4, step 2]{Bou_09a} 
 we obtained (\ref{faa}) for $f'\in C_b^{m-1}(\re)$, but here the assumption on $f$ implies only $f'\in C^{m-2}_b(\re)$.}
 We associate with $g$ the sequence $g_\nu:= \theta_\nu\ast g$, $\nu\geq 1$. It is easily seen that $(f\circ g_\nu)^{(\alpha)}\rightarrow ( f\circ g)^{(\alpha)}$ in the sense of distributions. Now we prove that the r.h.s. of (\ref{faa}) behaves similarly.\\
 
{\em Case 1: terms with $s<m$.}
Let us consider \be\label{lowterm}S(g):=(f^{(s)}\circ g)\, g^{(\gamma_{1})} \cdots
 g^{(\gamma_{s})}\,,\ee for a set of parameters satisfying (\ref{faa1}), and $s<m$.
 
 \begin{LEM}\label{low} The nonlinear operator $S$, given by (\ref{lowterm}), with $s<m$,  is well defined and continuous from $E$ to $L_p(\R)$.\end{LEM}

 {\bf Proof.} According to Gagliardo-Nirenberg (see e.g. \cite[(7), p.~6108]{Bou_09a}), we have the imbeddings
\[ \dot{A}^{m}_{p}(\re^{n}) \hookrightarrow \dot{W}^{N}_{pm/N}\,,\]
for all $N=1,\ldots , m-1$. Applying this to $N=|\gamma_r|$, for $r=1,\ldots, s$, and Proposition \ref{multili}, we obtain the 
continuity of the mapping $g\mapsto \, g^{(\gamma_{1})} \cdots
 g^{(\gamma_{s})}$ from $E$ to $L_p(\R)$. By Proposition \ref{alt}, we have $f^{(s)}\in C_b(\re)$. Hence we can apply Proposition \ref{gencont}, and conclude that $S$ is continuous. This ends up the proof of Lemma \ref{low}.\\
 
 As a consequence, it holds
 \[ \lim_{\nu\rightarrow \infty} (f^{(s)}\circ g_\nu)\, g_\nu^{(\gamma_{1})} \cdots
 g_\nu^{(\gamma_{s})} = (f^{(s)}\circ g)\, g^{(\gamma_{1})} \cdots
 g^{(\gamma_{s})}\,,\]
 in $L_p(\R)$.\\
 
 {\em Case 2 : terms with $s=m$.}
 Let us consider the operator
\[S(g):=(f^{(m)}\circ g)\,  \partial_{j_{1}}g\, \cdots\,\partial_{j_{m}}g\,,\]
with $j_{r}\in \{1,\ldots,n\}$ for $r=1,\ldots , m$.
For every subset $H$ of $\{1,\ldots, m\}$, let us denote
\[ D_H:= \sum_{k\in H} \partial_{j_k}\,.\]
According to Proposition \ref{algelem}, it holds
    \[S=\frac{(-1)^m}{m!}\sum_{H}(-1)^{|H| }\, S_{D_H}\,,\]
    see (\ref{S_D}) for the definition of $S_D$.
  By Lemma \ref{formula}, it holds $\lim_{\nu \rightarrow \infty} S(g_\nu)= S(g)$ in $L_p(\R)$. \\
  
  As a consequence of the above cases 1 and 2,  the formula (\ref{faa})
holds for all $g\in E$. This ends up the proof.

 \subsection{Proofs of Theorems \ref{contalg} and \ref{contw21} }

  {\em Proof of Theorem \ref{contalg}.}    We refer to the proof of Theorem \ref{compalg}. By Theorem \ref{contAF} and by the first imbedding given in (\ref{SAtoAF}), $T_{f\omega_r}$ is a continuous mapping
 from $\dot{W}^m_p\cap L_\infty (\R)$ to $\dot{W}^m_p(\R)$. Since $(f\omega_r)'\in L_\infty(\re)$, $T_{f\omega_r}$ is also  continuous 
 from $\dot{W}^m_p\cap L_\infty (\R)$ to  $L_\infty (\R)$. We conclude that $T_{f\omega_r}$ is continuous 
 from $\dot{W}^m_p\cap L_\infty (\R)$ to itself for every $r>0$. The continuity of $T_f$ follows by (\ref{inball}).\\

 {\em Proof of Theorem \ref{contw21}.} It  is similar to that of Theorem \ref{contAF}. We assume $n\not= 2$, since
the case $n=2$ is covered by Theorem \ref{contalg}.        
 By the estimations (\ref{estw21}), (\ref{n=1}), (\ref{esthomw21}),  and by the density of $\mathcal{D}(\re)$ into $L_1(\re)$, we can restrict ourselves to the case
     $f''\in  \mathcal{D}(\re)$. Notice that such a property implies $f'\in C^\infty_b(\re)$, hence
$\partial_j \partial_k(f\circ g) = U^{j,k}(g)+S^{j,k}(g)$, where $ U^{j,k}(g):= (f'\circ g)\partial_j\partial_k g$ and
      $S^{j,k}(g):= (f''\circ g)\,(\partial_jg)(\partial_kg)$, 
see \cite[step 2, p.~6111, and the proof of prop.~18, p.~6128]{Bou_09a}.
The continuity of $ U^{j,k}:E\rightarrow L_1(\R)$ follows by Proposition \ref{gencont}.
To prove the continuity of $S^{j,k}$, we introduce
\[ S_{D}(g):= (f''\circ g) \,(Dg)^2\,,\]
where $D$ is a first order differential operators with constant coefficients, and we set
\[ h(x):= \int_x^{+\infty} |f''(t)| \,\mathrm{d}t\,.\]
Then we must discuss according to $n$.\\

{\em Case $n>2$.} 
For all $g\in E$,  and $j=1,\ldots , n$, it holds $\partial_jg\in \dot{W}^1_1\cap L_r(\R)$, with $\frac{1}{r}:= 1-\frac{1}{n}$, see \cite[prop.~15]{Bou_09a}.
This implies the following property:
{\em for almost every $(x_1,\dots,x_{j-1}, x_{j+1},\ldots ,x_n)\in \re^{n-1}$, the function
\[t\mapsto \partial_jg(x_1,\dots,x_{j-1}, t,x_{j+1},\ldots ,x_n)\] is absolutely continuous on $\re$, with limit $0$ at $\pm\infty$.}
Thus we can make integrations by parts w.r.t. each coordinate, obtaining
\be\label{IP1} \|S_{D}(g)\|_1= \int_{\R} (h\circ g)\, (D^2g)\,\mathrm{d}x\,,\ee
for all $g\in E$. The continuity of $S^{j,k}:E\rightarrow L_1(\R)$ follows, exactly as in the proof Theorem \ref{contAF}.
By Sobolev's Theorem, the inhomogeneous space $W^2_1(\R)$ is imbedded into $E$. The continuity of
$T_f$ on $W^2_1(\R)$ follows at once if $f(0)=0$.\\

{\em Case $n=1$.} Now $D= {\mathrm d}/{\mathrm d}x$ and the formula (\ref{IP1}) becomes
\[\|S_{D}(g)\|_1= \int_{\re} (h\circ g)\, g''\,\mathrm{d}x\, - h\left(g(+\infty)\right)\,g'(+\infty)+  h\left(g(-\infty)\right)\,g'(-\infty)\,.\]
If $\lim g_\nu= g$ in $E$, then, by the proof of Lemma \ref{honestn=1}, $\lim g_\nu(x)= g(x)$ for every $x\in \re$, and also for $x=\pm\infty$. Since $h$ 
is continuous on $\re$, we conclude that $\lim  \|S_{D}(g_\nu)\|_1=  \|S_{D}(g)\|_1$. The rest of the proof is unchanged.

\section*{Conclusion}

Let us mention possible continuations of the present work:\\
1- {\em Generalization of theTheorem \ref{contm>0} to Sobolev spaces with fractional order of smoothness.}
The automatic continuity is known to hold in the following cases:
\begin{itemize}
\item Besov spaces $B^s_{p,q}(\R)$  with $0<s<1$, by interpolating between
$L_p$ and $W^1_p$, see \cite[5.5.2, thm.~3]{RS}.
\item Besov spaces $B^s_{p,q}(\re)$ and Lizorkin-Triebel spaces $F^s_{p,q}(\re)$
with $s>1+(1/p)$, $p,q\in [1,+\infty[$, see \cite[cor.~2]{BL} and \cite[thm.~8]{BS} (in those papers, they are some restrictions in case of Besov spaces, which have been removed
in \cite{114-1}).
 \end{itemize} 
The extension to the spaces on $\R$, for $n>1$ and $s>1$ noninteger, is completely open: the first difficulty is that we have not even a full characterization of functions which act by composition.

2- {\em Proof of the higher-order chain rule.}

In the proof of Theorem \ref{contAF}, we have established the formula (\ref{faa}) for all $g\in E$, but only for smooth functions $f$.
Could we generalize it to any $f$ s.t. $f'\in W^{m-1}_{L_{p,lu}}(\re)$? In this respect, we can refer to the partial
results of F.~Isaia \cite{IS}.

G\'erard Bourdaud
 
Universit\'e Paris Diderot, I.M.J. - P.R.G (UMR 7586)
       
Case 7012  
    
75205 Paris Cedex 13     \\                          
bourdaud@math.univ-paris-diderot.fr  \\

Madani Moussai

Laboratory of Functional  Analysis and  Geometry of  Spaces, 

M. Boudiaf University of M'Sila, 

28000 M'Sila, Algeria

mmoussai@yahoo.fr

\end{document}